\documentclass[12pt]{article}

\usepackage{psfig}
\usepackage{graphicx}
\usepackage{bm}
\usepackage{amsmath}
\usepackage{amssymb}
\usepackage[T2A]{fontenc}
\usepackage[cp866]{inputenc}

\newcommand{\e}{\varepsilon}
\newcommand{\pa}[2]{\frac{\partial #1}{\partial #2}}

\newtheorem{theorem}{Theorem}
\newtheorem{lem}{Lemma}
\newtheorem{rem}{Remark}
\newtheorem{cons}{Corollary}[theorem]

\makeatletter \@addtoreset{equation}{section}

\newcommand{\eq}[2]{\begin{equation}\label{#1}#2\end{equation}\addtocounter{syst}{1}}

\newcounter{syst}
\@addtoreset{syst}{section}
\renewcommand{\thesyst}{\thesection.\arabic{syst}}
\newcommand{\sy}{\refstepcounter{syst}\addtocounter{equation}{1}}

\makeatother \righthyphenmin=2

\author{Gorelyshev I.V.\qquad Neishtadt A.I.}
\title{On the adiabatic perturbation theory for systems with impacts}
\date{}

\begin{document}

\maketitle

{\small\it \centerline{Space Research Institute, Profsoyuznaya 84/32, Moscow 117997, Russia}

\vspace{2mm}

\centerline{E-mail: igor\_gor@iki.rssi.ru, aneishta@iki.rssi.ru} }

\vspace{1cm}

{\small{\bf Abstract.}We justify an applicability of the adiabatic
perturbation theory for three well known systems with impacts: a
ball between two slowly moving walls, a slowly irregular
waveguide, and an adiabatic piston.}

\vspace{5mm}

The adiabatic perturbation theory (see, for example,
\cite{arnold}, p. 200 -- 211) is used to describe the dynamics in
smooth Hamiltonian systems containing variables of two types: fast
and slow. The accuracy estimates of this theory are known. In a
number of cases, one can formally apply the procedure of this
theory to discontinuous Hamiltonian systems, in particular to
systems with impacts. However, the validity of such formal
approach does not follow from the theory for smooth systems. A
distinctive feature of considered problems is fast (instant)
variation of ``slow''{} variables at an impact. In the present
work we obtain accuracy estimates of the adiabatic perturbation
theory for systems with impacts. We consider three model problems
that belong to the three main classes of systems, where the
adiabatic perturbation theory is used: 1) the motion of a ball
between two slowly moving walls corresponds to the Hamiltonian
systems depending on slowly varying parameters; 2) the problem of
a slowly irregular waveguide corresponds to the Hamiltonian
systems slowly depending on some of variables ; 3) the problem of
an adiabatic piston corresponds to the Hamiltonian systems with
slow and fast motions. We discuss these problems in order to
demonstrate the effectiveness of the approach based on the
adiabatic perturbation theory: the dynamics in these problems is
well known, and it was described earlier by other methods (in
particular, in \cite{sinai} for problem 3)).

The first approximation of the procedure of the adiabatic perturbation theory leads to a conclusion
that the system has an adiabatic invariant (an approximate integral). Such a conclusion is often
used for systems with impacts, but its validity is derived from direct calculations (see \cite{A},
Sec. 52). Higher approximations of the procedure of the adiabatic perturbation theory for systems
with impacts are formally considered in \cite{neishtadt}. The perturbation theory for non-smooth
Hamiltonian systems is discussed in works \cite{markeev1,markeev2}, yet only in the case of
continuous phase variables. In systems with impacts some of phase variables are discontinuous at a
point of reflection.

Dealing with systems with impacts, one usually considers Poincare map in order to obtain results
analogous to those obtained for smooth systems by means of higher approximations of the adiabatic
perturbation theory. If this map is smooth, one can either use the procedure of the adiabatic
perturbation theory for smooth maps (see, for example, \cite{J}), or use an artificial approach,
namely, to describe the map as a Poincare map for an auxiliary smooth Hamiltonian system and apply
the procedure of the adiabatic perturbation theory to this system \cite{N}.

An advantage of the approach described in the present paper is
that one can deal with Hamiltonians of original systems, and
consider systems with impacts similarly to smooth systems. It is
no need to construct and deal with Poincare maps (that can turn
out to be non-smooth as, for example, in the problem of an
adiabatic piston, sec.4). Under this approach calculations become
simpler and more concise.

\section{Fermi-Ulam model}

The problem of oscillations of a particle between two slowly
moving walls is called the Fermi-Ulam problem (or model)
(Fig.\ref{fermi})\cite{Z}. Let the left wall be fixed (this
assumption is made to shorten calculations) and the right wall
slowly change its position. The distance between the walls is
$d(\e t)\ge \mbox{const}> 0$, where $t$ is the time variable, $\e
> 0$ is a small parameter. Introduce slow time variable $\tau$, $\tau=\e
t$. Let function $d(\cdot)\in C^{\infty}$. We consider the problem
on a slow time interval that either does not depend on $\e$ or
grows as $\e$ decrease. In the second case we shall suppose that
function $d$ and its derivatives are bounded on the real axis. Let
the mass of the particle be equal to $1$, then the velocity of the
particle is equal to its momentum. One can describe the dynamics
of the particle in the following way: between the walls the
particle has a constant velocity $v$; the velocity of the particle
change sign after an impact with the left wall, and transforms as
$v_1=2\dot{d}-v$ after an impact with the right moving wall. The
Hamiltonian of the particle is well defined between walls and
describes free motion of the particle with a constant momentum
$v$.

One can interpret the problem as the motion of a particle in a
potential that equals $0$ between the walls and $+\infty$ in the
other space. The right wall is moving slowly. Therefore, it is
useful to consider the problem when the right wall is fixed
($d=\mbox{const}$). The phase portrait of this system is shown in
Fig.\ref{phase}. On the portrait one can define ``action-angle''{}
variables $(I,\phi)$ in a standard way. ``Action''{} is the area
bounded by a phase trajectory divided by $2\pi$:
$I=\frac{1}{2\pi}\cdot |v|\cdot 2d=\frac{d}{\pi}\sqrt{2E}$, where
$E$ is the Hamiltonian of the particle. Hence,
$E=\frac{\pi^2I^2}{2d^2}$. The angle (phase) $\phi$ is the
uniformly changing angle variable on the trajectory;
$\phi=2\pi\frac{t}{T}$, where $t$ is the time interval
corresponding to the motion of the particle from the initial point
to a given point and $T$ is the period of motion.

Suppose the phase is zero at the left wall; then

\sy
$$\phi =
\left\{
\begin{array}{lll}
\pi\frac{x}{d},\qquad \phi\in (0,\,\pi),\\
\pi(2-\frac{x}{d}),\qquad \phi\in (\pi,\, 2\pi),\\
\end{array}
\right.\eqno (\thesyst)
$$
\label{sys1.1}where $x$ is the distance between the particle and
the left wall.

If the phase is defined by (\ref{sys1.1}), then generating function $W=W(x,I,d)$ of the form

\sy
$$
W= \left\{\begin{array}{lll}
\pi\frac{Ix}{d},\qquad \phi\in (0,\,\pi),\\
\pi(2-\frac{Ix}{d}),\qquad \phi\in (\pi,\, 2\pi)\\
\end{array}\right.\eqno (\thesyst)
$$transforms variables $(v,x)$ to the ``action-angle''{} variables. \label{sys1.2}

Use generating function $W$ to make a change of variables in the
real system with right wall moving slowly. The motion of the
particle between the walls is described by the Hamiltonian system
with Hamiltonian $H=E+\pa{W}{t}$:

\eq{1.1}{H=\frac{\pi^2I^2}{2d^2}-\frac{I\dot{d}}{d}f(\phi),}

where

\eq{1.2}{f(\phi)=\left\{\begin{array}{lll}
\phi,\qquad \phi\in (0,\,\pi),\\
\phi-2\pi,\qquad \phi\in (\pi ,\, 2\pi).\\
\end{array}\right.}

Hamiltonian (\ref{1.1}) completely describes the motion of a
particle between the walls.

At an impact with the fixed wall, the value of the ``action''{}
variable is preserved. The variation of the ``action''{} variable
after an impact with the moving wall is

\eq{1.4}{I_{+}-I_{-}=-\frac{2d\dot{d}}{\pi},}where $I_-$ and $I_+$
are the values of $I$ prior to and after the impact
correspondingly.

Let us describe an impact by means of Hamiltonian (\ref{1.1}). An
interaction with the wall is instantaneous. Hence, it is likely
that the actual value of the variation of the ``action''{}
variable after an impact with the moving wall can be calculated in
the following way: fix the moment of time at the impact, consider
Hamiltonian (\ref{1.1}) as a Hamiltonian of a one d.o.f.
Hamiltonian system, and obtain the required value of the
``action''{} variable using the energy conservation law.

\begin{lem}
This approach properly describes jump of the ``action''{} variable
at an impact with the moving wall. \label{l1.1}
\end{lem}

{\bf Proof.}

Suppose the value of time variable is fixed. Set the value of
energy prior to an impact equal to the value of energy after an
impact:

$$\frac{\pi^2 I_{-}^2}{2d^2}-\frac{\pi I_{-}\dot d}{d}=\frac{\pi^2 I_{+}^2}{2d^2}+\frac{\pi I_{+}\dot d}{d}.$$
Hence, we obtain

$$\frac{\pi (I_{+}^2-I_{-}^2)}{2d}=-\dot d (I_{+}+I_{-}),$$
from which immediately follows (\ref{1.4}).

As applied to the considering problem, the basic statement of the adiabatic perturbation theory can
be formulated as follows.

\begin{theorem}
After canonical change of variables $(I,\phi)\,\mapsto\, (\hat I,\hat\phi)$ with generating
function

\eq{1.w}{W=\hat I\phi+\e S(\hat I,\phi,\tau,\e),\qquad S = S_1+\e S_2+\ldots+\e^{r-2}S_{r-1},}
where $r$ is any fixed natural number, Hamiltonian (\ref{1.1}) takes the form

\eq{1.h}{{\cal H}={\cal H}_{\Sigma, r
}(\hat{I},\tau,\e)+\e^{r}H_{r}(\hat{I},\phi(\hat{I},\hat\phi,\tau,\e),\tau,\e),}
$${\cal H}_{\Sigma, r}=E+\e {\cal H}_1+\ldots+\e^{r-1}{\cal H}_{r-1}.$$

Here $S_i=S_i(\hat I,\phi,\tau)$, ${\cal H}_i={\cal H}_i(\hat I,\tau)$. Functions $S_i$, ${\cal
H}_i$ are $C^{\infty}$ in $\hat I,\tau$; functions $S_i$ are continuous with respect to $\phi$.

\label{t1}
\end{theorem}

{\bf Proof.}

We use the standard procedure of the adiabatic perturbation theory (see, for example,
\cite{arnold}). Make a canonical change of variables with generating function (\ref{1.w}) in the
system with Hamiltonian (\ref{1.1}). The variables transform according to the following expressions

\eq{1.f}{I=\hat I+\e\pa{S}{\phi},\qquad
\hat\phi=\phi+\e\pa{S}{\hat{I}}.}

New Hamiltonian ${\cal H}$ has the form

\eq{1.nh}{{\cal H}=H(I,\phi,\tau,\e)+\e^2\pa{S(\hat
I,\phi,\tau,\e)}{\tau}.}

Define functions $S_i$ so that the new Hamiltonian has the form (\ref{1.h}). Substitute expressions
(\ref{1.f}) in (\ref{1.nh}) and equal terms of the same order in $\e$. Thus, we find the set of
equations from which $S_i$,${\cal H}_i$ can be found. For example, in the first order in $\e$ we
get

$${\cal H}_1(\hat I,\tau)=\pa{E}{I}\pa{S_1}{\phi}-\frac{Id'}{d}f(\phi).$$

A prime denotes derivative with respect to $\tau$. After averaging over $\phi$ one obtains ${\cal
H}_1=0$. Therefore, $S_1$ can be found as a quadrature, which turns out to be continuous with
respect to$\phi$. One can choose $S_1$ so that its average is $0$. Similarly, one can derive
functions $S_i$ that are continuous with respect to $\phi$ and $C^{\infty}$ with respect to $\hat
I$,$\tau$. The required properties of functions ${\cal H}_i$,$H_{r}$ follow.

\begin{cons}
The value of variable $\hat{I}$ is preserved along a solution with
an accuracy $O(\e^{r})$ on time interval $O(\e^{-k})$ for any
prefixed natural number $k$. The value of variable $\hat{\phi}$ is
defined by an integrable Hamiltonian system with Hamiltonian
${\cal H}_{\Sigma, r}$ with an accuracy $O(\e^{r-k})$ on the same
time interval. \label{c1}
\end{cons}

{\bf Proof.}

Consider the change of variables $(I,\phi)\,\mapsto\,(\tilde
I,\tilde\phi)$ generated by function $\tilde W=\tilde
I\phi+\e\tilde{ S}(\tilde I,\phi,\tau,\e)$, $\tilde S=S_1+\e
S_2+\ldots+\e^{r+k-2}S_{r-1}$.

According to Theorem \ref{t1}, the Hamiltonian function in the new variables has the form:

\eq{1.7}{\tilde{\cal H}={\cal H}_{\Sigma, r+k}(\tilde
I,\tau,\e)+\e^{r+k}H_{r+k}(\tilde I,\phi,\tau,\e).}

The expressions for the change of variables give $\hat I-\tilde I
= O(\e^{r})$. Hence, below we consider variation of variables
$\tilde I$ and $\hat\phi$. Variation of value $\tilde I$ between
the impacts with the walls is $O(\e^{r})$. Consider variation of
value $\tilde I$ at an impact with the moving wall. Lemma
\ref{l1.1} describes the way to calculate the value of variable
$I$ after an impact with the wall. In order to determine variation
of variable $\tilde{I}$ one substitutes
$I=I(\tilde{I},\phi,\tau,\e)$ in Hamiltonian $H$ in (\ref{1.1})
and notes that $H={\cal H}-\e^2\pa{\tilde{S}}{\tau}$. Thus, one
obtains the following expression

$${\cal H}_{\Sigma, r+k}(\tilde{I}_+,\tau,\e)+\e^{r+k}H_{r+k}(\tilde{I}_+,\pi+0,\tau,\e)-\e^2\pa{\tilde{S}(\tilde
I_+,\pi,\tau,\e)}{\tau} ={\cal H}_{\Sigma, r+k
}(\tilde{I}_-,\tau,\e)+$$
$$+\e^{r+k}H_{r+k}(\tilde{I}_-,\pi-0,\tau,\e)-\e^2\pa{\tilde{S}(\tilde
I_-,\pi,\tau,\e)}{\tau}.$$

Here $\tilde I_-$ and $\tilde I_+$ are the values of $\tilde I$
prior to and after the impact respectively. We note that the
function $\tilde S$ is continuous with respect to $\phi$. Hence,
we find
$$\left(\pa{E(\tilde I_-,\tau)}{\tilde
I}+O(\e^2)\right)(\tilde I_+-\tilde I_-)=O(\e^{r+k}).$$

Therefore $\tilde I_+-\tilde I_-=O(\e^{r+k})$.

The particle collides with the walls not more than $O(\e^{k})$
times. Therefore the value of $\tilde I$ varies as $O(\e^r)$.
Hence the part of the statement concerning variable $\tilde I$ is
proved.

Consider now variation of $\hat\phi$. On any time interval between
the impacts with the moving wall value of $\hat\phi$ deviates from
the solution of equation $\dot{\hat{\phi}}=\pa{{\cal H}_{\Sigma,
r}}{\hat I}$, $\hat I = \mbox{const}$ by $O(\e^{r})$. At a
collision with the moving wall value of variable $\hat\phi$ jumps.
Value of the jump can be found from the formulas for change of
variables (\ref{1.f}) and depends on variation of $\hat I$.
Functions $S_m$ are continuous with respect to $\phi$. Hence, the
jump in the phase at an impact with the wall is $O(\e^{r+1})$.
This completes the proof of the statement.

\begin{cons}
Value of variable $I$ is preserved with an accuracy $O(\e)$ on a
time interval $O(\e^{-k})$ for any prefixed natural number $k$.
\label{c2}
\end{cons}

\begin{cons}
The formulas for change of variables $\hat I,\hat\phi \,\mapsto\,
I,\phi $ and approximate expressions $\hat I = \mbox{const}$,
$\dot{\hat{\phi}}=\pa{{\cal H}_{\Sigma, r}}{\hat I}$ describe the
behaviour of variables $I,\phi$ with an accuracy $O(\e^r)$ for $I$
and $O(\e^{r-k})$  for $\phi$ on the time interval $O(\e^{-k})$.
\label{c3}
\end{cons}

\begin{rem}
Lemma \ref{l1.1}, theorem \ref{t1} and corollaries \ref{c1} - \ref{c3} can be used in case of
motion in potential $U(x,\tau)\in C^{\infty}$ between the walls, and the particle in adiabatic
approximation collides with the walls at a nonzero velocity.
\end{rem}

\begin{rem}
In \cite{KT} systems with impacts are investigated in the
following way. Instead of system with impacts, a smooth system
with a repulsive potential strong in a thin layer (the thickness
of the layer $\delta << 1$, the gradient of the potential $\sim
1/{\delta}$) is considered. Properties of the system with impacts
are derived from that of the smooth system proceeding to limit as
$\delta\to 0$. Following this approach one could define such a
potential for Fermi-Ulam problem (and for problems considered
below, sec. 2,3), use the perturbation theory for the smooth
system and obtain corresponding accuracy estimates for the system
with impacts. However, one should consider uniformity of accuracy
estimates in $\delta$. The method defined above allows to avoid
this problem.
\end{rem}

\section{A slowly irregular planar waveguide}

In this section we consider slowly irregular planar waveguides.
The problem is to find the trace for a ray of light propagating in
a slowly irregular waveguide with reflective walls. A slow
irregularity means that the width of the waveguide is changing
slowly along the waveguide's length \cite{KO}.

\begin{figure}
\begin{tabular}{p{8cm}p{8cm}}
\psfig{file=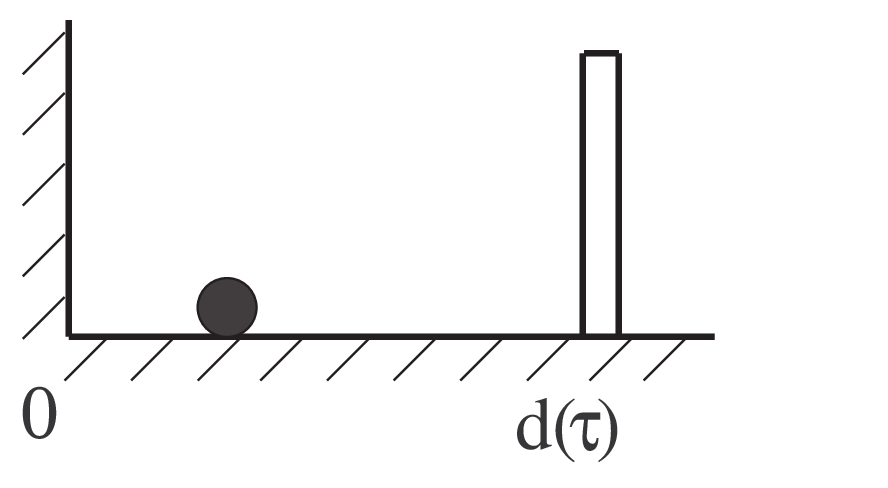,width=200pt,angle=0} \caption{\small{A ball
between slowly moving walls.}}
\label{fermi} & \psfig{file=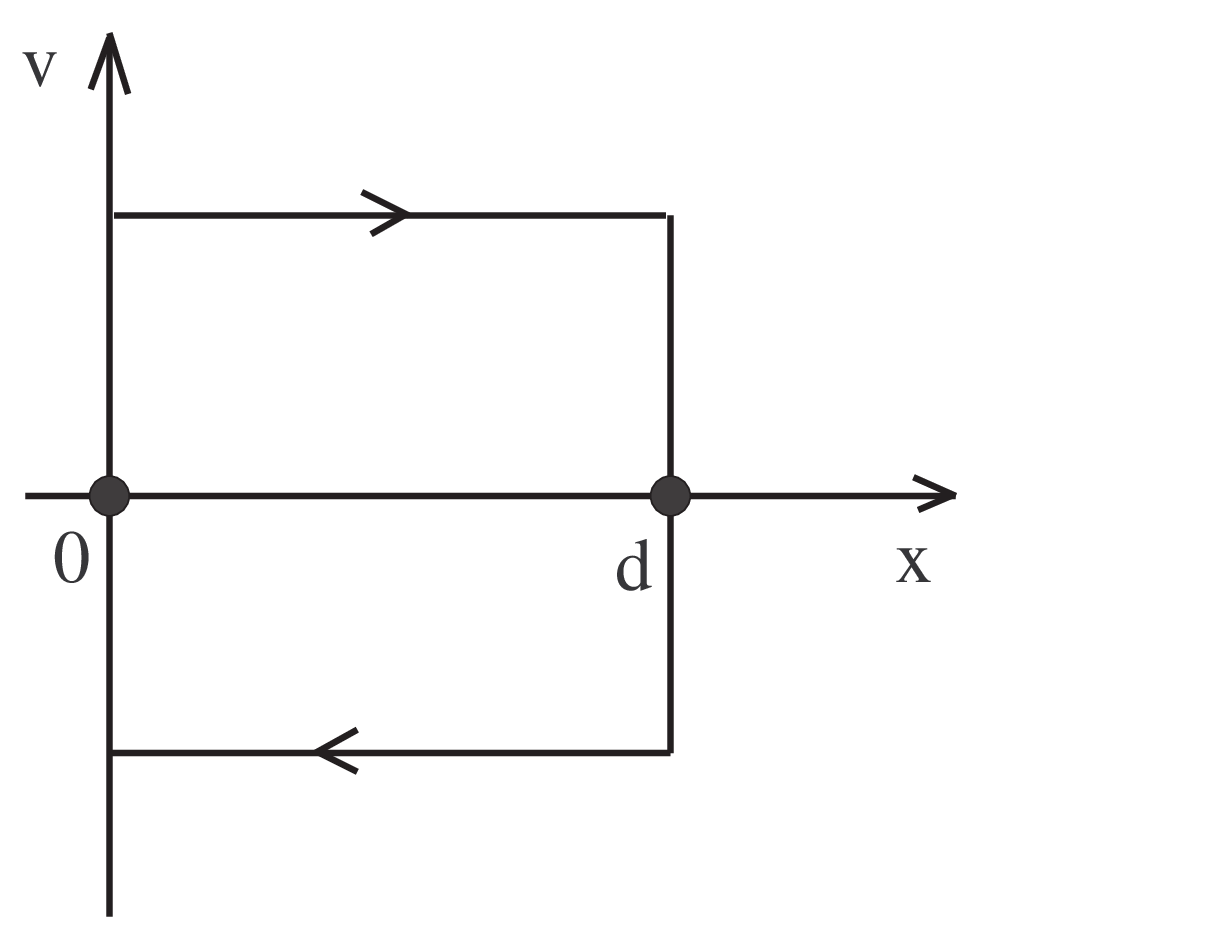,width=200pt,angle=0} \caption{\small{A phase portrait for the motion of the particle.}} \label{phase}\\
\psfig{file=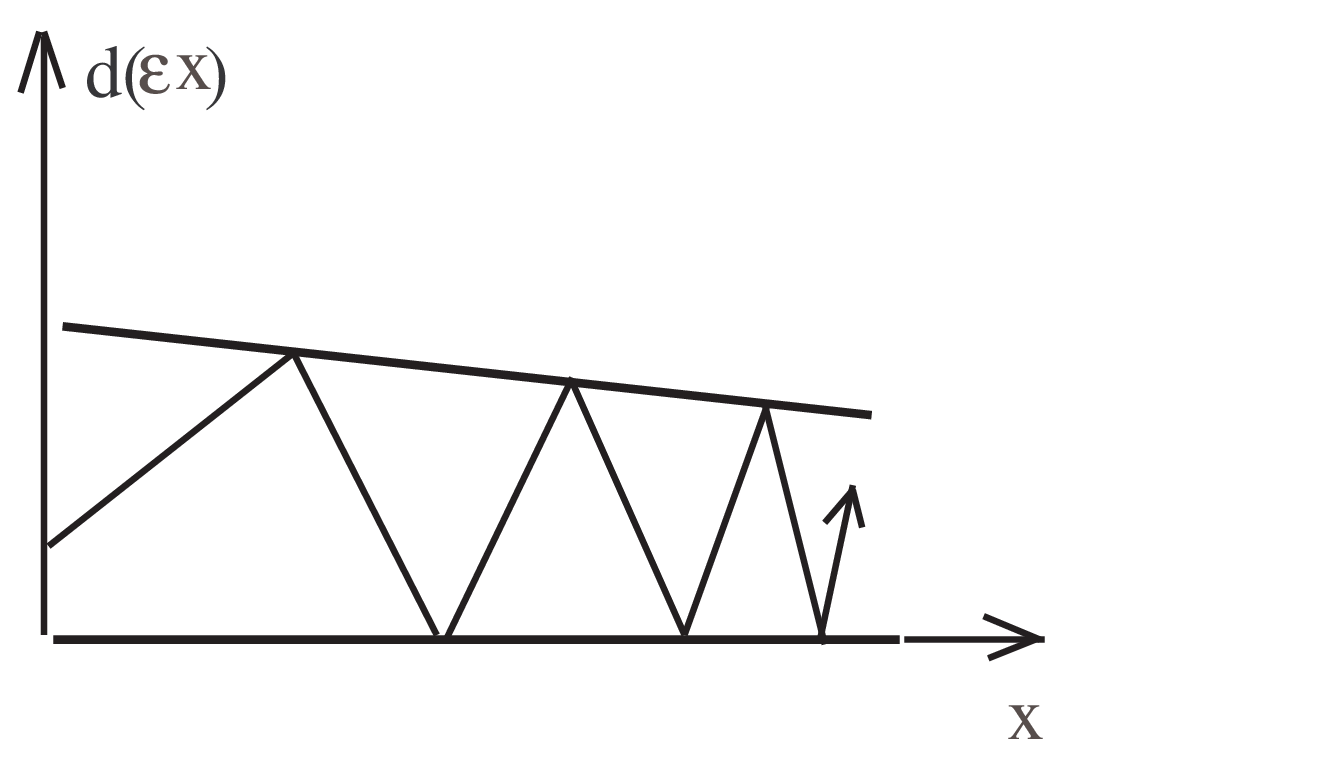,width=200pt,angle=0} \caption{\small{A planar
slowly irregular waveguide.}} \label{wave} &
\psfig{file=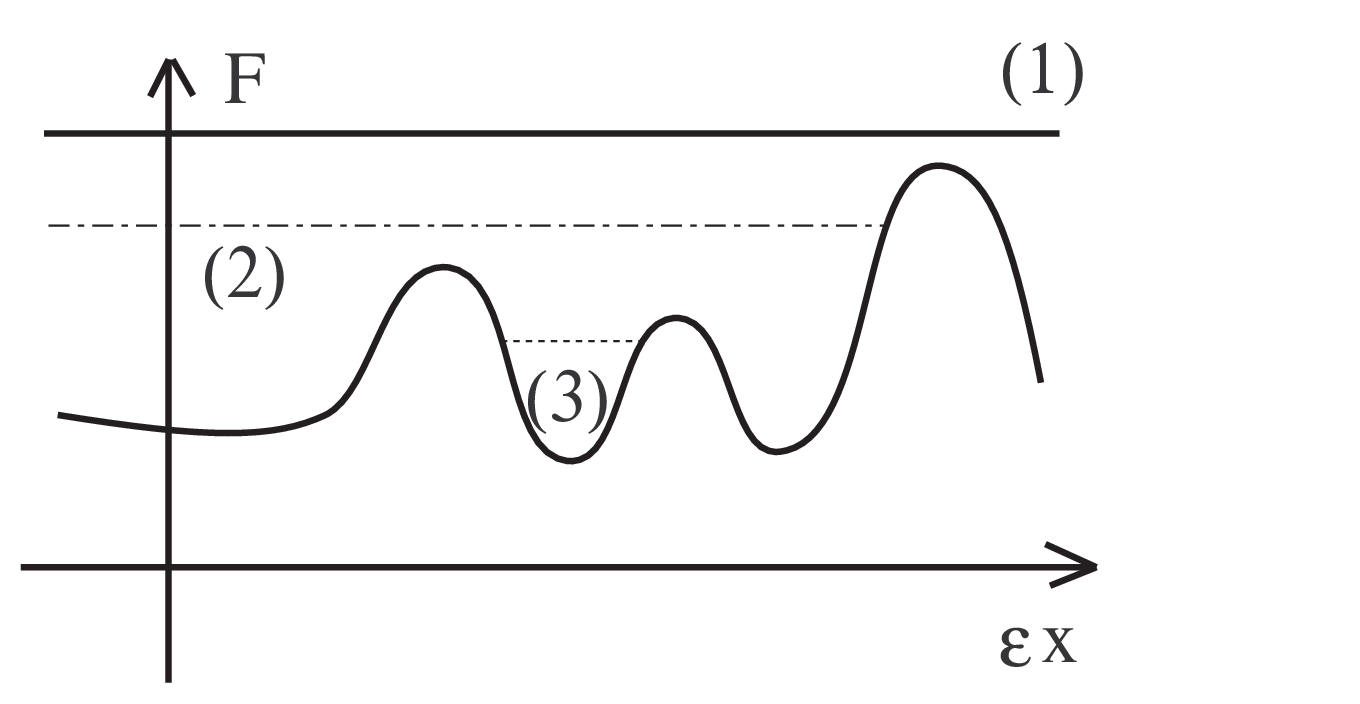,width=200pt,angle=0} \caption{\small{The
plot of function $2U$
and levels of function $F$.}} \label{guide}\\
\psfig{file=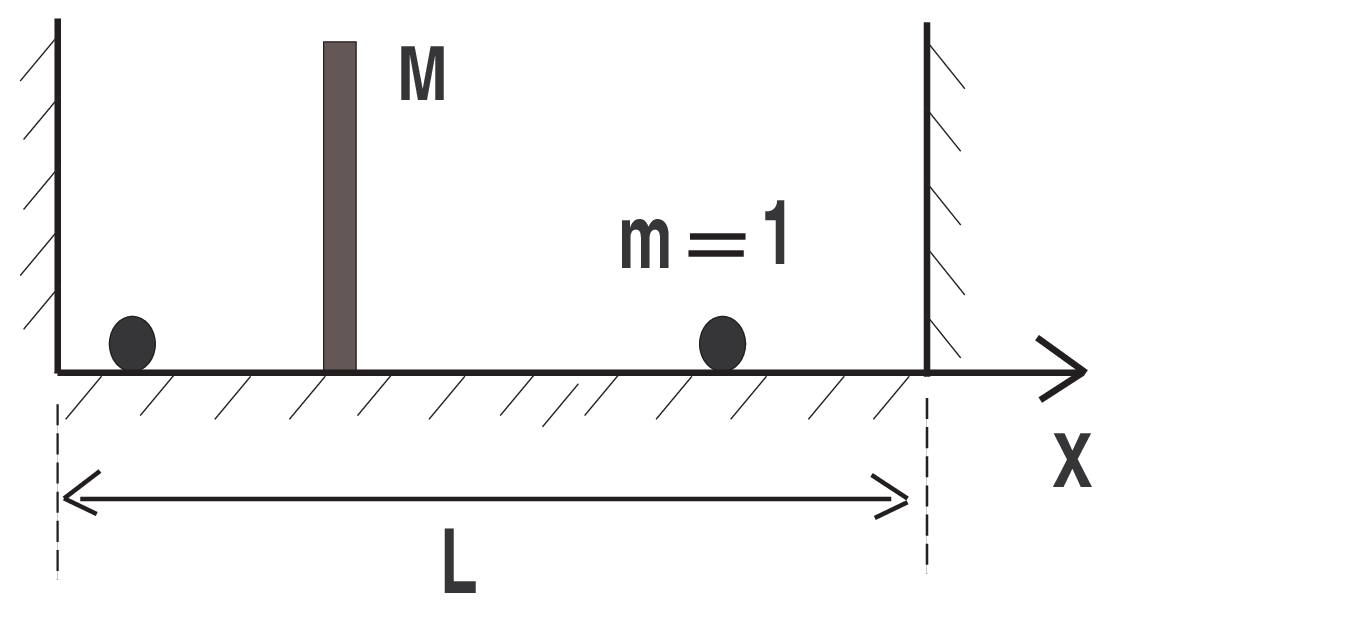,width=200pt,angle=0} \caption{\small{The
system of a piston and two particles.}}
\label{wall} & \psfig{file=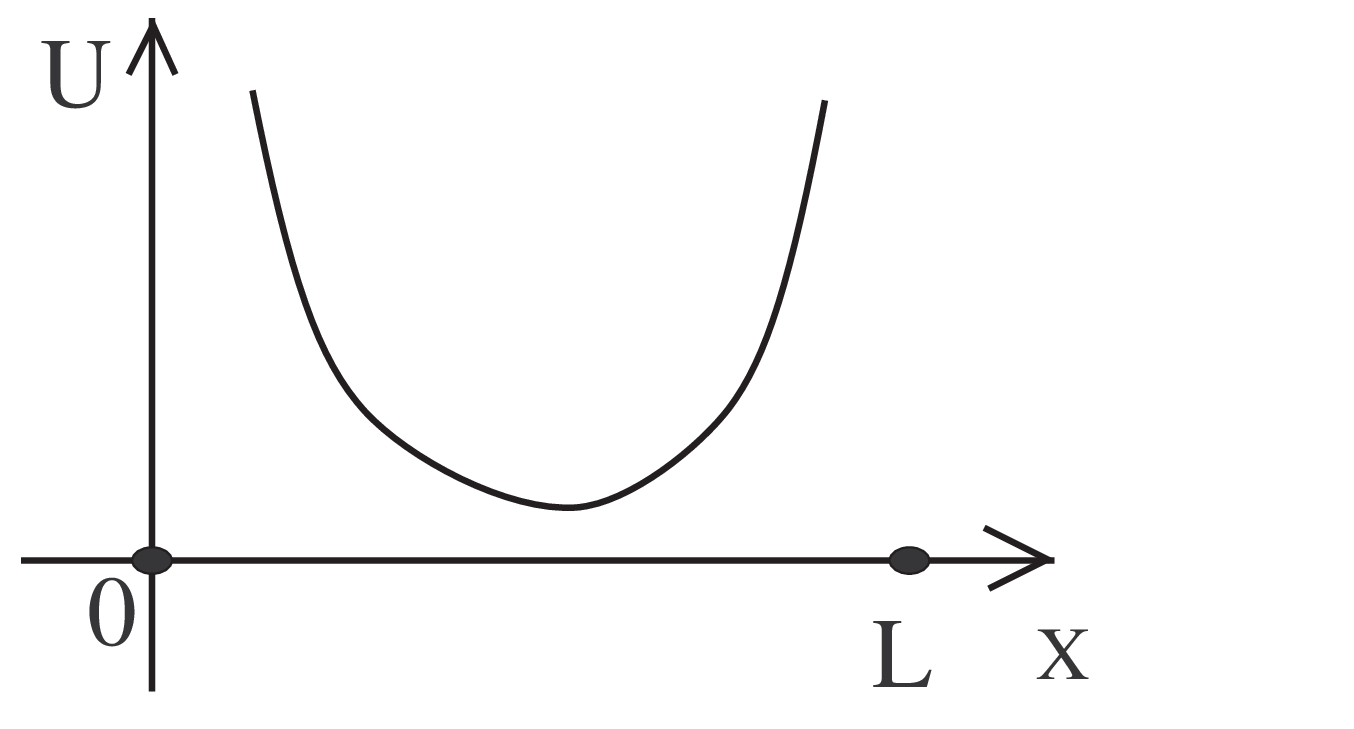,width=200pt,angle=0} \caption{\small{The effective potential energy of the piston.}} \label{potential}\\
\end{tabular}
\end{figure}

Assume, for simplicity, that one wall of the waveguide coincides
with $x$ axis. The position of another wall in the plane $Oxy$ is
defined by formula $y=d(X)$, where $X=\e x$, $\e$ is a small
parameter (Fig.\ref{wave}). Let $d(\cdot)\in C^{\infty}$. Let us
consider the problem on a time interval $O(\e^{-k})$, where $k$ is
any prefixed natural number. Assume that function $d$, its
derivatives of any order and function $1/d$ are uniformly bounded
on the real axis.

One can describe the propagation of rays of light in a medium
using the Hamiltonian system with Hamiltonian

\eq{3.1}{H=p_x^2+p_y^2-n^2(x,y),} where one should consider only
zero energy level $H=0$ \cite{KO}. Here $p_{x,y}$ are the
variables canonically conjugated to the coordinates $x,y$ and
$n(x,y)$ is the refraction index. In our case inside the waveguide
$n=1$. The Hamiltonian system is not well defined at the
reflective walls and one can use only conservation laws there. In
the unperturbed system $X=\mbox{const}$ and the projection of a
phase trajectory on $p_y,\, y$ plane looks similar to the phase
trajectory of a particle in Fig.\ref{phase}. Like in the previous
section, one can make a canonical change of variables $(p_y,\, y)$
to ``action-angle''{} variables $(I,\,\phi)$: $I=\frac{|p_y|\cdot
d(X)}{\pi}$; if phase $\phi$ is zero at the bottom wall, then

\sy
$$\phi =
\left\{
\begin{array}{lll}
\pi\frac{y}{d(X)},\qquad \phi\in (0,\,\pi),\\
\pi(2-\frac{y}{d(X)}),\qquad \phi\in (\pi,\, 2\pi).\\
\end{array}
\right.\eqno (\thesyst)
$$
\label{sys3.1} Using standard definition of the phase
(\ref{sys3.1}) one obtains the generating function for the change
of variables $(p_y,\, y)$ to the ``action-angle''{} variables in
the following form

\sy
$$
W(I,y,X)= \left\{\begin{array}{lll}
\pi\frac{Iy}{d},\qquad \phi\in (0,\,\pi),\\
\pi(2-\frac{Iy}{d}),\qquad \phi\in (\pi,\, 2\pi).\\
\end{array}\right.\eqno (\thesyst)
$$
\label{sys3.2}

Make a canonical change of variables ($p_x$, $x$, $p_y$, $y$)
$\mapsto$ ($\hat{p}_x$, $x$, $I$, $\phi$) with generating function
$S=\hat{p}_x x+W(I,y,\e x)$ in the exact (perturbed) system. The
Hamiltonian in the new variables is

\eq{3.1a}{H=\frac{\pi^2I^2}{d^2}+\left(\hat{p}_x-\e\frac{I\,
d'}{d}f(\phi)\right)^2-1,} where $f(\phi)$ is defined in
(\ref{1.2}), the prime denotes a derivative with respect to $X$.

At a reflection at the bottom wall value of the ``action''{} is
preserved. At a reflection at the upper wall it transforms
according to the following expression

\eq{3.2}{I_1=I-\frac{2\e d'}{1+\e^2{d'}^2}(\e d'I+\frac{d
p_x}{\pi}).}

\begin{rem}
Value $\hat{p}_x$ is preserved at a reflection. Formula
(\ref{3.2}) defines the law of transformation for the angle of the
ray of light arriving to the bottom wall. This angle transforms
after an impact with the upper wall as
$\alpha_1=\alpha-2\,arctg(\e d')$, thus implying (\ref{3.2}).
\end{rem}

\begin{lem}
Consider Hamiltonian (\ref{3.1a}) and fix variables
$(\hat{p}_x,\,x)$. Value of the Hamiltonian is preserved when
value of $\phi$ passes through $\pi$. This law of conservation
defines the actual value of ``action''{} variable $I$ after
reflection at the upper wall of the waveguide. \label{l3}
\end{lem}

{\bf Proof.}

Set equal values of Hamiltonian (\ref{3.1a}) before and after
reflection of the ray at the upper wall:

$$\frac{\pi^2I^2}{d^2}+\left(\hat{p}_x-\e\frac{I\, d'}{d}\pi\right)^2-1=
\frac{\pi^2{I_1}^2}{d^2}+\left(\hat{p}_x+\e\frac{I_1\, d'}{d}\pi\right)^2-1.$$After transformations
obtain

$$\frac{\pi^2I^2}{d^2}(1+\e^2{d'}^2)-2\e\hat{p}_x\frac{I\, d'}{d}\pi=
\frac{\pi^2{I_1}^2}{d^2}(1+\e^2{d'}^2)+2\e\hat{p}_x\frac{I_1\,
d'}{d}\pi.$$As value of the ``action''{} variable is positive, we
find

$$I_1=I-\frac{2\e d'}{1+\e^2{d'}^2}\frac{\hat{p}_x d}{\pi}.$$Substituting the expression for momentum $p_x=\hat{p}_x-\e\frac{I\, d'}{d}\pi$, we finally get formula

$$I_1=I-\frac{2\e d'}{1+\e^2{d'}^2}(\e d'I+\frac{dp_x}{\pi}),$$
coinciding with (\ref{3.2}). The lemma is proved.

Let us formulate the basic statement of the adiabatic perturbation
theory for the considered problem.

\begin{theorem}

After canonical change of variables $(\hat{p}_x, x, I, \phi)$
$\mapsto$ $(\tilde{p}_x,\tilde x, \tilde I,\tilde \phi)$ with
generating function

\eq{3.3}{W=\tilde I\phi+\tilde{p}_xx+\e S(\tilde
I,\phi,\tilde{p}_x,\e x,\e),\qquad S=S_1+\e S_2+\ldots+\e^{r-2}
S_{r-2},} where $r$ is any prefixed natural number, Hamiltonian
(\ref{3.1a}) takes the following form

\eq{2.H}{H={\cal H}_{\Sigma, r}(\tilde I,\tilde{p}_x,\e\tilde
x,\e)+\e^{r}H_{r}(\tilde I,\phi,\tilde{p}_x,\e x,\e),}

$${\cal H}_{\Sigma, r}=\frac{\pi^2\tilde{I}^2}{d^2}+\tilde{p}_x^2-1+\e{\cal H}_1+\ldots+\e^{r-1}{\cal H}_{r-1}.$$
Here functions $S_i$, ${\cal H}_i$ are $C^{\infty}$ with respect
to $\tilde I,\tilde{p}_x,\e\tilde x$; $S_i$ are continuous
functions of $\phi$.

\label{t3}
\end{theorem}

{\bf Proof.}

Make a canonical change of variables $(I,\phi,\hat{p}_x,x)\,\mapsto\, (\tilde
I,\tilde\phi,\tilde{p}_{x},\tilde x)$ with generating functions (\ref{3.3}). The variables are
transformed as:

\sy
$$
\left\{\begin{array}{lll}
I=\tilde I+\e\pa{S}{\phi},\\
\tilde \phi=\phi+\e\pa{S}{I},\\
\hat{p}_x=\tilde{p}_{x}+\e^2\pa{S}{\e x},\\
\e\tilde x=\e x+\e^2\pa{S}{\tilde{p}_{x}}.\\
\end{array}
\right. \eqno (\thesyst)
$$
\label{sys3.2a}

Following the standard procedure of the adiabatic perturbation theory, define functions $S_m$ such
that the Hamiltonian in the new variables does not depend on $\phi$ up to terms of order $O(\e^m)$
and $S_m$ are continuous and periodic with respect to $\phi$. The procedure is as follows:
substitute formulas for the change of the ``action''{} and longitudinal variables ($\hat{p}_x$,$x$)
from (\ref{sys3.2a}) to Hamiltonian (\ref{3.1a}) and choose the generating function so that the
Hamiltonian in the new variables has form (\ref{2.H}).

Thus, one can recurrently define functions $S_m$. For example, an
equation for $S_1$ looks as follows:

$$\pa{S_1}{\phi}=\frac{\tilde{p}_x d(X)d'}{\pi^2}f(\phi).$$

One can choose functions $S_m$ in such a way that they have zero average with respect to $\phi$.
Finally, the Hamiltonian takes the form we need. Functions $S_i,H_i$ obviously have all the
properties mentioned in the theorem.

\vskip 0.3cm

Consider the solutions of Hamiltonian equations with Hamiltonians
(\ref{3.4}) and ${\cal H}_{\Sigma, r}(J,p,\e x,\e)$ corresponding
to initial conditions $\tilde{I}(0)=J(0)
,\,\tilde\phi(0)=\psi(0),\,\tilde{p}_x(0)=p(0),\,
\tilde{x}(0)=x(0)$, where $\psi$ is the phase variable conjugated
to variable $J$.

\begin{cons}

Value of variable $\tilde I$ is preserved with an accuracy
$O(\e^r)$ on a time interval $\e^{-k}$ for any prefixed natural
number $k$. The projection of the ray's trajectory on
$\tilde{p}_x,\e\tilde x$ plane on such an interval lies in the
$O(\e^r)$-neighbourhood of curve ${\cal H}_{\Sigma, r}(J,p ,\e x,
\e)= 0$. Behaviour of variables $\e \tilde
x,\tilde{p}_x,\tilde\phi$ is described by the solution of the
system with Hamiltonian ${\cal H}_{\Sigma, r}$ with accuracies
$O(\e^{r-k+1}), O(\e^{r-k+1})$, and  $O(\e^{r-k})$ respectively
(under some natural additional conditions given
below).\label{c2.1}
\end{cons}

{\bf Proof.}

Similarly to the previous section, consider an auxiliary change of
variables with generating function

$$ W=\bar I\phi+\bar{p}_xx+\e S_1(\bar
I,\phi,\bar{p}_x,\e x)+\ldots+\e^{r+k-1} S_{r+k-1}(\bar
I,\phi,\bar{p}_x,\e x).
$$
According to Theorem \ref{t3}, Hamiltonian (\ref{3.1a}) takes the
form:
$$H={\cal H}_{\Sigma, r+k}(\bar
I,\bar{p}_x,\e\bar x,\e)+\e^{r+k}H_{r+k}(\bar I,\phi,\bar{p}_x,\e
x,\e),$$ where $\phi$ and $x$ are considered as functions of new
variables $(\bar I,\bar\phi,\bar{p}_x,\bar x)$. Formulas for the
change of variables give $\tilde I -\bar I = O(\e^r)$, $\tilde\phi
- \bar\phi = O(\e^r)$, $\tilde x-\bar x = O(\e^r)$ and
$\tilde{p}_x - \bar{p}_x = O(\e^{r+1})$. Therefore, it is
sufficient to consider behaviour of variables $\bar
I,\tilde\phi,\tilde{p}_x,\tilde x$.

Consider an approximate conservation of variable $\bar I$. The
ray's Hamiltonian in the new variables is

\eq{3.4}{H={\cal H}_{\Sigma, r+k}(\bar I,\bar{p}_x,\e\bar
x,\e)+\e^{r+k}H_{r+k}(\bar I,\phi,\bar{p}_x,\e x,\e).} On each
time interval between the reflections of the ray at the upper wall
variable $\bar I$ vary according to Hamiltonian equation

$$\dot{\bar{I}}=-\e^{r+k}\pa{H_{r+k}}{\bar{\phi}}.$$ After summing up over all such intervals, it
gives an accuracy of conservation $O(\e^r)$.

It follows from (\ref{3.4}) that values of variables
$(\bar{p}_x,\,\bar x)$ jump at a reflection of the ray at the
upper wall. However, the generating function for the change of
variables is continuous with respect to $\phi$. Therefore,
formulas for the change of variables give
$\Delta\bar{p}_x=O(\e^2\Delta\bar{I})$,
$\Delta\e\bar{x}=O(\e^2\Delta\bar{I})$, where $\Delta\bar{I}$,
$\Delta\bar{p}_x$, $\Delta\e\bar{x}$ are jumps of variables $\bar
I$, $\bar{p}_x$, $\e\bar x$ after reflection of the ray at the
upper wall. Using Lemma \ref{l3}, one estimates value of the jump
of variable $\bar{I}$ as $O(\e^{r+k})$. The ray reflects at the
upper wall not more than $O(\e^{-k})$ times on the time interval
$\e^{-k}$. Therefore, total variation of variable $\bar I$ is
$O(\e^r)$. Thus, the part of the theorem concerning variable $\bar
I$ is proved.

Projection of the rays trajectory on $\tilde{p}_x,\e\tilde x$
plane is in $O(\e^r)$-neighbourhood of the curve ${\cal
H}_{\Sigma, r}(J,p ,\e x, \e)= 0$ because along the trajectory
${\cal H}_{\Sigma, r}(J,\tilde{p}_x,\e\tilde x,\e)= O(\e^r)$.

Values of variables $p,x,\psi$ change continuously, yet values of variables
$\bar{p}_x,\bar{x},\bar\phi$ jump at a reflection of the ray at the upper wall. It was noted above
that values of jumps of variables $\bar{p}_x,\bar{x},\bar \phi$ are determined by the value of jump
of $\bar{I}$ at the reflection of the ray at the upper wall as
$\Delta\bar{p}=O(\e^2\Delta\bar{I})=O(\e^{r+k+2})$ $\Delta\bar{x}=O(\e\Delta\bar{I})=O(\e^{r+k+1})$
and $\Delta\bar\phi=O(\e\Delta\bar{I})=O(\e^{r+k+1})$. Thus, total variations of values of
variables ($\bar p_x$,$\bar x$,$\bar\phi$) due to reflections are $O(\e^{r+2})$ for variable
$\bar{p}_x$ and $O(\e^{r+1})$ for variables $\bar{x}$ и $\bar\phi$. Total variations of variables
$\tilde p_x$,$\tilde x$,$\tilde \phi$ due to reflections are the same.

Consider the zero approximation of Hamiltonian ${\cal H}_{\Sigma,
r}:$
$${\cal H}_{\Sigma, r}^0\stackrel{def}{=}F=\frac{\pi^2 I^2}{d^2(\e x)}+p^2$$
The Hamiltonian system with Hamiltonian $\frac{1}{2}{\cal
H}_{\Sigma, r}^0$ describes the motion of a particle in potential
$U=\frac{\pi^2 I^2}{2 d^2(\e x)}$. Assume (see Fig.\ref{guide})
that for the ray under consideration value of function $F$
satisfies one of the following conditions: 1) it is larger than
value of any local maximum of $2U$ or 2) the ray can reflect at a
potential hump only once or 3) the ray moves between two potential
humps (resonator).

Consider case 1). The motion of the ray has a given direction,
variables $\tilde x$ and $x$ vary monotonically. One can consider
variables $\tilde p_x$ and $p$ as functions of variables $\e\tilde
x$ and $\e x$ respectively. At $x=\tilde x$ values of variables
$\tilde p_x$ and $p$ differ by $O(\e^r)$. It takes different time
for variables $\tilde x$ and $x$ to reach a given value $x_*$.
This difference can be estimated as follows:

$$
\int\limits_{x(0)}^{x_* }\left(\frac{1}{\dot{\tilde
{x}}}-\frac{1}{\dot x}\right)dx =O(\e^{r-k}).
$$
Hence, for a given moment of time $t$ one obtains: $\e \tilde
x(t)-\e x(t) = O(\e^{r-k+1}), \, \tilde p_x(t)- p(t) =
O(\e^{r-k+1})$.

The deviation of phases $\tilde\phi-\psi$ at $x=\tilde x=x_*$ is :
\eq{3.6}{\tilde\phi-\psi=\int\limits_{x(0)}^{x_*
}\left(\frac{\omega(\tilde I ,\tilde {p}_x,\e
x,\e)}{\dot{\tilde{x}}}-\frac{\omega(J,p,\e x,\e)}{\dot
x}\right)dx+O(\e^{r-k}) =O(\e^{r-k}) ,}where we introduced
$\omega=\pa{{\cal H}_{\Sigma, r}}{I}$. Therefore, for a given
moment of time $t$ one gets: $\tilde\phi(t)-\psi(t)=O(\e^{r-k})$.

In case 2) the ray can change its direction of propagation only
once. Here one can consider two different regions and use either
variable $x$ or $p$ as monotonically changing variable in the
corresponding region. In case 3) the waveguide is configured as a
resonator. One can introduce an ``angle'' {} variable in the phase
portrait of the system with Hamiltonian ${\cal H}_{\Sigma,
r}(I,p,x)$ as a monotonically changing variable. All arguments
here are the same as in case 1).

\begin{cons}
Value of variable $I$ is preserved with an accuracy $O(\e)$ on a
time interval $\e^{-r}$ for any prefixed natural number $r$.
\label{c2.2}
\end{cons}

\begin{cons}
Behaviour of variables $I,\phi,\hat{p}_x,x$ can be found from the change of variables $\tilde
I,\tilde\phi,\tilde{p}_x,\tilde{x} \,\mapsto\, I,\phi,\hat{p}_x,x$ together with approximate
expressions defining motion in the system with Hamiltonian ${\cal H}_{\Sigma, r}(\tilde
I,\tilde{p}_x, \e\tilde{x})$. An accuracy of such a description is the same as the accuracy for
variables $\tilde I,\tilde\phi,\tilde{p}_x,\tilde{x}$ according to corollary \ref{c2.1}.
\label{c2.3}
\end{cons}

\begin{rem}
The Hamiltonian system with Hamiltonian

\eq{3.5}{H=\frac{\pi^2I^2}{d^2}+\hat{p}_x^2-1}considered on the energy level $H=0$ describes the
ray's trajectory with an accuracy $O(\e)$ on time intervals of order $\e^{-1}$.

\end{rem}

\begin{rem}

Suppose that there is nonhomogeneous medium between the walls such
that a refraction index $n(\e x,y)\in C^{\infty}$ and trace of the
ray in adiabatic approximation is inclined to the walls at
reflection at a nonzero angle; then Lemma \ref{l3}, Theorem
\ref{t3} and its corollaries \ref{c2.1} -- \ref{c2.3} are still
valid.
\end{rem}

\section{Dynamics of a massive piston in a gas of light particles}

The problem of an adiabatic piston is an important model in statistical mechanics. It is considered
in the context of attempts to derive thermodynamics laws from the laws of mechanics (see, for
example, \cite{sinai},\cite{lieb}). The corresponding system consists of a container with a massive
cylindrical piston and a gas of identical light particles that move independently elastically
colliding with the walls of the container and with the piston.

Let the mass of a particle be equal to $1$. Length of the
container subtracting thickness of the piston $L$ and the number
of gas particles are of the order $1$. Mass of the piston $M$ is
large in comparison with mass of the gas. Suppose the piston is at
rest at the initial moment of time; then energy of the system does
not depend on mass of the piston. Hence, one can estimate energy
of the piston as $O(1)$ and its typical velocity as
$O(\frac{1}{\sqrt{M}})$. Therefore, it is useful to introduce a
small parameter $\e=\frac{1}{\sqrt{M}}$. Let us discuss the
problem on a time interval of order $\e^{-1}$.

Below, to shorten the calculations, we consider the case when there is only one particle on each
side of the piston (Fig.\ref{wall}). One can easily generalize following arguments for the case of
any prefixed number of particles. At the end of the section we give the result for the general
case. Without loss of generality, assume that the particles' velocities are parallel to the axis of
the container; thus, the motion is one-dimentional. Let indexes $l,r$ correspond to variables of
left and right particles respectively. The variables without indexes correspond to the piston.
Thus, full energy of the system is of the form:

\eq{1}{E=\e^2\frac{P^2}{2}+\frac{p_l^2}{2}+\frac{p_r^2}{2},}where
$p_{l,r}$ are momenta of particles, $P$ is momentum of the piston.
Denote as $x_{l,r}$ distances between particles and the left wall
of the container, as $X$ distance between the piston and the left
wall of the container. The distance for the right particle is
calculated subtracting the thickness of the piston.

Suppose the piston is fixed. For each particle make a canonical
change of variables from $(p,x)$ to ``action-angle''{} variables
$(I,\phi)$ the same way as in Fermi-Ulam problem (Sec.$1$). Let
phases $\phi_{l,r}$ be zero at the walls of the container; then

\sy
$$\phi_l =
\left\{
\begin{array}{lll}
\pi\frac{x_l}{X},\qquad \phi\in (0,\,\pi),\\
\pi(2-\frac{x_l}{X}),\qquad \phi\in (\pi,\, 2\pi),\\
\end{array}
\right. \qquad
  \phi_r =
\left\{
\begin{array}{lll}
\pi\frac{L-x_r}{L-X},\qquad \phi\in (0,\,\pi),\\
\pi(2-\frac{L-x_r}{L-X}),\qquad \phi\in (\pi,\, 2\pi).\\
\end{array}
\right. \eqno (\thesyst)
$$
\label{sys3.p}

Conjugated ``action''{} variables are $I_l=\frac{|p_l|X}{\pi}$,
$I_r=\frac{|p_r|(L-X)}{\pi}$. If phases are defined by
(\ref{sys3.p}), generating functions for the change of variables
$(p,x)$ to ``action-angle''{} have the following forms:

$$
S_l(I_l,x_l,X)= \left\{\begin{array}{lll}
\pi\frac{I_lx_l}{X},\qquad \phi_l\in (0,\,\pi),\\
\pi I_l(2-\frac{x_l}{X}),\qquad \phi_l\in (\pi,\, 2\pi),\\
\end{array}\right.
$$
$$
S_r(I_r,x_r,X)= \left\{\begin{array}{lll}
\pi\frac{I_r(L-x_r)}{L-X},\qquad \phi_r\in (0,\,\pi),\\
\pi I_r(2-\frac{L-x_r}{L-X}),\qquad \phi_r\in (\pi,\, 2\pi).\\
\end{array}\right.
$$

Make the change of variables with generating function $W=\hat{P}X+S_l(I_l,x_l,X)+S_r(I_r,x_r,X)$ in
the exact system, when the piston is moving. The Hamiltonian in the new variables has the form:

\eq{3}{{\cal H} =
\e^2\frac{1}{2}\left(\hat{P}-\frac{I_l}{X}f(\phi_l)+\frac{I_r}{L-X}f(\phi_r)\right)^2+\frac{\pi^2
I_l^2}{2X^2}+\frac{\pi^2 I_r^2}{2(L-X)^2},}where function
$f(\phi)$ defined by (\ref{1.2}).

Let us find jump of the ``action''{} variable at an impact. It is
natural to fix values of all variables except the variables of the
colliding particle, obtain the Hamiltonian for the particle and
consider the conservation law for the value of the Hamiltonian.
Let us formulate the corresponding lemma for the left particle
(the lemma for the right particle is formulated analogously). Let
the prime denote value of the variable after an impact.

\begin{lem} This approach defines actual value of the ``action''{} variable after an impact:

\eq{8}{I'_l=I_l-\frac{2PX}{\pi M}+\frac{2}{M+1}(\frac{P X}{\pi
M}-I_l).} \label{l2}
\end{lem}

\begin{rem}
Velocities $v$ and $V$ of the particles of mass $1$ and $M$
transform after an elastic collision as follows:

$$
\left\{\begin{array}{lll}
V'=V+\frac{2}{M+1}(v-V),\\
v'=2V-v+\frac{2}{M+1}(v-V).
\end{array}
\right.
$$

``Action''{} of the particle and momentum of the piston are
$I_l=\frac{|v|X}{\pi}$, $P=MV$, from which follows (\ref{8}).
\end{rem}

\begin{rem}
If two particles collide with the piston at the same moment, the
further dynamics is not defined. Measure of initial conditions
corresponding to such collisions equals $0$. Therefore, we do not
consider these initial condition.
\end{rem}

{\bf Proof.}

An impact is the the passage of $\phi_l$ through the value $\pi$. Equating values of the energy
prior to- and after an impact, we find

\eq{5}{\frac{\pi^2I_l^2}{2X^2}\frac{M+1}{M}-\frac{\pi I_l}{M
X}(\hat{P}+\frac{I_rf(\phi_r)}{L-X})=
\frac{\pi^2I_l'^2}{2X^2}\frac{M+1}{M}+\frac{\pi I'_l}{M
X}(\hat{P}+\frac{I_rf(\phi_r)}{L-X}).}

After calculations one obtains value of the ``action''{} variable:

\eq{6}{I'_l=I_l-\frac{2}{M+1}\frac{X}{\pi}(\hat{P}+\frac{I_rf(\phi_r)}{L-X}).}

Let $P$ denote value of momentum of the piston prior to an impact;
then

\eq{7}{I'_l=I_l-\frac{2}{M+1}\frac{X}{\pi}(P+\frac{\pi I_l}{X}).}

Finally

$$I'_l=I_l-\frac{2PX}{\pi M}+\frac{2}{M+1}(\frac{P X}{\pi
M}-I_l),$$

and the lemma is proved.

\begin{theorem}
Values of variables $I_{l,r}$ are preserved with an accuracy
$O(\e)$ on a time interval $\e^{-1}$.
\end{theorem}

{\bf Proof.}

Introduce the normalized momentum of the piston by expression
$\check P=\e\hat P$ and consider Hamiltonian (\ref{3}) with an
accuracy  $O(\e^2)$:

\eq{9}{{\cal H} = \frac{\check{P}^2}{2}+\frac{\pi^2
I_l^2}{2X^2}+\frac{\pi^2 I_r^2}{2(L-X)^2}-\e
\check{P}\frac{I_l}{X}f(\phi_l)+\e
\check{P}\frac{I_r}{L-X}f(\phi_r)+O(\e^2).}

Make a canonical change of variables ($I_{l,r}$, $\phi_{l,r}$,
$\e^{-1}\check{P}$, $X$) $\mapsto$ ($\tilde{I}_{l,r}$,
$\tilde{\phi}_{l,r}$, $\e^{-1}\tilde{P}$, $\tilde X$), eliminating
dependence of the Hamiltonian on the phase with an accuracy
$\e^2$. The generating function for this change of variables has
the form:

\eq{10}{W = \frac{1}{\e}\tilde{P}X+\tilde I_l\phi_l+\tilde
I_r\phi_r+ \e S_l(\tilde I_l,\phi_l,\tilde{P},X)+ \e S_r(\tilde
I_r,\phi_r,\tilde{P},X).}

Consider the expressions connecting old and new variables:

\sy
$$
\left\{\begin{array}{lll}
I_{l,r}=\tilde I_{l,r}+\e\pa{S_{l,r}}{\phi_{l,r}},\\
\tilde \phi_{l,r}=\phi_{l,r}+\e\pa{S_{l,r}}{I_{l,r}},\\
P=\tilde{P}+\e^2\pa{S_l}{X}+\e^2\pa{S_r}{X},\\
\tilde X=X+\e^2\pa{S_l}{\tilde{P}}+\e^2\pa{S_r}{\tilde{P}}.\\
\end{array}
\right. \eqno (\thesyst)
$$
\label{sys2.2}

Substitute these expressions in Hamiltonian (\ref{9}) and choose
functions $S_{l,r}$ so that the Hamiltonian in the new variables
does not depend on the phases in the first order in $\e$.

Therefore, functions $S_{l,r}$ satisfy equations

\eq{11}{\pa{S_{l}}{\phi_{l}}=\frac{\tilde{P}X
f(\phi_l)}{\pi^2},\qquad \pa{S_r}{\phi_r}=-\frac{\tilde{P}(L-X)
f(\phi_r)}{\pi^2}.}

It follows from (\ref{11}) that functions $S_{l,r}$ are defined up
to an arbitrary function of variables $\tilde{I}_{l,r}, \tilde{P},
X$. Let this function be equal to $0$ and call variables
$\tilde{I}_{l,r}$ the ``improved actions''{}. It is easy to see
from (\ref{sys2.2}), (\ref{11}), that the new phase coincides with
the old one ($\tilde\phi_{l,r}=\phi_{l,r}$), and the Hamiltonian
has a discontinuity at $\tilde{\phi}_{l,r}=\pi$.

Lemma \ref{l2} and formulas (\ref{sys2.2}), (\ref{11}) allow to
find variation of the ``improved action''{} after an impact as
$O(\e^2)$. Variation of the ``improved actions''{} between
collisions is also $O(\e^2)$. One can estimate the number of
collisions on a time interval $\e^{-1}$ as $O(\e^{-1})$. Hence, on
a time interval $\e^{-1}$ ``improved actions''{} vary by $O(\e)$.
The ``action''{} variables are related to ``improved actions''{}
by formulas for the change of variables
$I_{l,r}=\tilde{I}_{l,r}+\e\pa{S_{l,r}}{\phi_{l,r}}$. Therefore,
``action''{} variables vary by $O(\e)$ on the same time interval
and the theorem is proved.

\begin{cons}Let
\eq{12}{H=\e^2\frac{
P^2}{2}+\frac{\pi^2I_l^2}{2X^2}+\frac{\pi^2I_r^2}{2(L-X)^2},\qquad
I_{l,r}=\mbox{const},}where $I_{l,r}$ are initial values of the
``action'' variables.  The Hamiltonian system with Hamiltonian
(\ref{12}) describe the behaviour of variables $\e P, X$ with an
accuracy $O(\e)$ on a time interval $\e^{-1}$.

\end{cons}

One can similarly consider the system with any prefixed number of
particles. Value of the ``action''{} variable of each particle is
preserved with an accuracy $O(\e)$ on a time interval $\e^{-1}$.
The Hamiltonian system with Hamiltonian

\eq{13}{H=\e^2\frac{ P^2}{2M}+\frac{\pi^2}{2X^2}I_{\Sigma, l
}^2+\frac{\pi^2}{2(L-X)^2}I_{\Sigma, r}^2,}describes behaviour of
variables $\e P,X$ with an accuracy $O(\e)$ on the same time
interval. Here $I_{\Sigma, r}^2$ and $I_{\Sigma, l}^2$ are the
sums of squared initial values of ``action''{} variables for the
left particles and for the right particles respectively. (This
result was first derived in \cite{sinai} by a different method.)

\begin{rem}
In the considered approximation the piston oscillates in potential
$U=\frac{\pi^2}{2\tilde{X}^2}I_{\Sigma, l
}^2+\frac{\pi^2}{2(L-\tilde{X})^2}I_{\Sigma, r}^2$
(Fig.\ref{potential}).

\end{rem}

\section{Acknowledgements}

The work was partially supported by RFBR (03-01-00158,
НШ136.2003.1)

\newpage

\end{document}